\documentclass[5p]{elsarticle} % 
	\usepackage{amssymb}
	\usepackage{amsmath}
	\usepackage{graphicx}	
	\usepackage{tikz,pgfplots}	
	\pgfplotsset{compat=newest}
	\usepackage{pgf}
	\usepackage{graphicx}
	\usepackage{lineno}		
	\usepackage{draftwatermark}
	\SetWatermarkScale{5}

\newtheorem{definition}{Definition}[section]

\newtheorem{proposition}{Proposition}[section]

\begin{document}
\begin{frontmatter}
\title{Uniqueness of the Power Flow Solutions in Low Voltage Direct Current Grids}
\author{Alejandro Garces}%
\address{email: alejandro.garces@utp.edu.co \\ Universidad Tecnologica de Pereira. \\ 
AA: 97 - Post Code: 660003 - Pereira, Colombia}

\begin{abstract}

Power flow in a low voltage direct current grid (LVDC) is a non-linear problem just as its counterpart ac. This paper demonstrates that, unlike in ac grids, convergence and uniqueness of the solution can be guaranteed in this type of grids. The result is not a linearization nor an approximation, but an analysis of the set of non-linear algebraic equations, which is valid for any LVDC grid regardless its size, topology or load condition. Computer simulation corroborate the theoretical analysis.

\end{abstract}
		

\end{frontmatter}
\section{Introduction}

\subsection{Motivation}

Low voltage direct current (LVDC) is a promising technology for urban distribution systems, micro-grids, data centers, traction power systems and shipboard power systems \cite{overview}. It presents advantages in terms of reliability, efficiency, controllabiliy, power density and loadability \cite{LVDC_distribution,local_dc_distribution}.

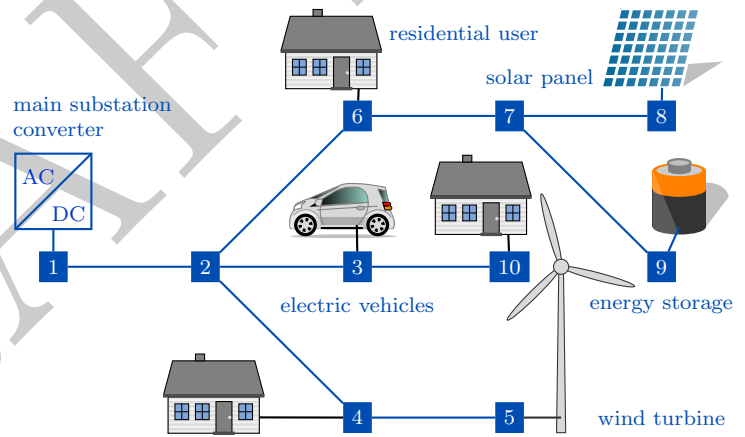
\begin{figure}[tb]
\centering
\footnotesize
\begin{tikzpicture}[x=1.0mm, y = 1.0mm,blue!70!green]
\node at (0,0) [draw,thick,fill, text =white] (N1) {1};
\node at (20,0) [draw,thick,fill, text =white] (N2) {2};
\node at (40,0) [draw,thick,fill, text =white] (N3) {3};
\node at (40,-20) [draw,thick,fill, text =white] (N4) {4};
\node at (60,-20) [draw,thick,fill, text =white] (N5) {5};
\node at (40,20) [draw,thick,fill, text =white] (N6) {6};
\node at (60,20) [draw,thick,fill, text =white] (N7) {7};
\node at (80,20) [draw,thick,fill, text =white] (N8) {8};
\node at (80,0) [draw,thick,fill, text =white] (N9) {9};
\node at (60,0) [draw,thick,fill, text =white] (N10) {10};
\draw[-,thick] (N1) -- (N2);
\draw[-,thick] (N2) -- (N3);
\draw[-,thick] (N3) -- (N10);
\draw[-,thick] (N2) -- (N6);
\draw[-,thick] (N2) -- (N4);
\draw[-,thick] (N6) -- (N7);
\draw[-,thick] (N7) -- (N9);
\draw[-,thick] (N7) -- (N8);
\draw[-,thick] (N4) -- (N5);

\node at (64,25) {solar panel};
\node at (54,31) {residential user};
\node at (40,-5) {electric vehicles};
\node at (80,-20) {wind turbine};
\node at (80,-5) {energy storage};
% Usuario residencial 1
\begin{scope}[gray, xshift=104, yshift=80, scale=0.3]
\draw[black,  fill=gray] (-22,2) -- +(4,15) -- +(40,15) -- +(44,0) -- cycle;
\draw[black,  fill=gray] (-8,17) -- +(0,3) -- +(4,3) -- +(4,0) -- cycle;
\draw[gray!30,top color=gray!30, bottom color=white] (-20,1) rectangle +(40,-2);
\draw[gray!30,top color=gray!30, bottom color=white] (-20,-1) rectangle +(40,-2);
\draw[gray!30,top color=gray!30, bottom color=white] (-20,-3) rectangle +(40,-2);
\draw[gray!30,top color=gray!30, bottom color=white] (-20,-5) rectangle +(40,-2);
\draw[gray!30,top color=gray!30, bottom color=white] (-20,-7) rectangle +(40,-2);
\draw[gray!30,top color=gray!30, bottom color=white] (-20,-9) rectangle +(40,-2);
\draw[gray!30,top color=gray!30, bottom color=white] (-20,-11) rectangle +(40,-2);
\draw[gray!30,top color=gray!30, bottom color=white] (-20,-13) rectangle +(40,-2);
\draw[black, fill=blue!70!green!50] (-17,-9) rectangle +(6,7);
\draw[black] (-14,-9) -- +(0,7);
\draw[black] (-17,-5.5) -- +(6,0);
\draw[black, fill=blue!70!green!50] (-8,-9) rectangle +(6,7);
\draw[black] (-5,-9) -- +(0,7);
\draw[black] (-8,-5.5) -- +(6,0);
\draw[black,fill=blue!70!green!50] (11,-9) rectangle +(6,7);
\draw[black] (14,-9) -- +(0,7);
\draw[black] (11,-5.5) -- +(6,0);
\draw[black, fill=gray] (1,-15) rectangle +(7,14);
\draw[gray!30,fill] (7,-8) circle (0.5);
\draw[black] (-20,1) rectangle +(40,-16);
\draw[black!70, fill] (-22,1) rectangle +(44,1);
\draw[black!70, fill] (0,-15.5) rectangle +(9,1);
\draw[-,thick, black] (12,-15) -- (N6);
\end{scope}

% Usuario residencial 2
\begin{scope}[gray, xshift=160, yshift=25, scale=0.3]
\draw[black,  fill=gray] (-22,2) -- +(4,15) -- +(40,15) -- +(44,0) -- cycle;
\draw[black,  fill=gray] (-8,17) -- +(0,3) -- +(4,3) -- +(4,0) -- cycle;
\draw[gray!30,top color=gray!30, bottom color=white] (-20,1) rectangle +(40,-2);
\draw[gray!30,top color=gray!30, bottom color=white] (-20,-1) rectangle +(40,-2);
\draw[gray!30,top color=gray!30, bottom color=white] (-20,-3) rectangle +(40,-2);
\draw[gray!30,top color=gray!30, bottom color=white] (-20,-5) rectangle +(40,-2);
\draw[gray!30,top color=gray!30, bottom color=white] (-20,-7) rectangle +(40,-2);
\draw[gray!30,top color=gray!30, bottom color=white] (-20,-9) rectangle +(40,-2);
\draw[gray!30,top color=gray!30, bottom color=white] (-20,-11) rectangle +(40,-2);
\draw[gray!30,top color=gray!30, bottom color=white] (-20,-13) rectangle +(40,-2);
\draw[black, fill=blue!70!green!50] (-17,-9) rectangle +(6,7);
\draw[black] (-14,-9) -- +(0,7);
\draw[black] (-17,-5.5) -- +(6,0);
\draw[black, fill=blue!70!green!50] (-8,-9) rectangle +(6,7);
\draw[black] (-5,-9) -- +(0,7);
\draw[black] (-8,-5.5) -- +(6,0);
\draw[black,fill=blue!70!green!50] (11,-9) rectangle +(6,7);
\draw[black] (14,-9) -- +(0,7);
\draw[black] (11,-5.5) -- +(6,0);
\draw[black, fill=gray] (1,-15) rectangle +(7,14);
\draw[gray!30,fill] (7,-8) circle (0.5);
\draw[black] (-20,1) rectangle +(40,-16);
\draw[black!70, fill] (-22,1) rectangle +(44,1);
\draw[black!70, fill] (0,-15.5) rectangle +(9,1);
\draw[-,thick, black] (12,-15) -- (N10);
\end{scope}

% Usuario residencial 3
\begin{scope}[gray, xshift=60, yshift=-50, scale=0.3]
\draw[black,  fill=gray] (-22,2) -- +(4,15) -- +(40,15) -- +(44,0) -- cycle;
\draw[black,  fill=gray] (-8,17) -- +(0,3) -- +(4,3) -- +(4,0) -- cycle;
\draw[gray!30,top color=gray!30, bottom color=white] (-20,1) rectangle +(40,-2);
\draw[gray!30,top color=gray!30, bottom color=white] (-20,-1) rectangle +(40,-2);
\draw[gray!30,top color=gray!30, bottom color=white] (-20,-3) rectangle +(40,-2);
\draw[gray!30,top color=gray!30, bottom color=white] (-20,-5) rectangle +(40,-2);
\draw[gray!30,top color=gray!30, bottom color=white] (-20,-7) rectangle +(40,-2);
\draw[gray!30,top color=gray!30, bottom color=white] (-20,-9) rectangle +(40,-2);
\draw[gray!30,top color=gray!30, bottom color=white] (-20,-11) rectangle +(40,-2);
\draw[gray!30,top color=gray!30, bottom color=white] (-20,-13) rectangle +(40,-2);
\draw[black, fill=blue!70!green!50] (-17,-9) rectangle +(6,7);
\draw[black] (-14,-9) -- +(0,7);
\draw[black] (-17,-5.5) -- +(6,0);
\draw[black, fill=blue!70!green!50] (-8,-9) rectangle +(6,7);
\draw[black] (-5,-9) -- +(0,7);
\draw[black] (-8,-5.5) -- +(6,0);
\draw[black,fill=blue!70!green!50] (11,-9) rectangle +(6,7);
\draw[black] (14,-9) -- +(0,7);
\draw[black] (11,-5.5) -- +(6,0);
\draw[black, fill=gray] (1,-15) rectangle +(7,14);
\draw[gray!30,fill] (7,-8) circle (0.5);
\draw[black] (-20,1) rectangle +(40,-16);
\draw[black!70, fill] (-22,1) rectangle +(44,1);
\draw[black!70, fill] (0,-15.5) rectangle +(9,1);
\draw[-,thick, black] (20,-8) -- (N4);
\end{scope}
% Solar panel
\begin{scope}[ very thick, scale=0.5, xshift=410, yshift=135]	
	\fill[gray!50] (0,0) -- +(5,5) -- +(32,7) -- +(21,0) -- cycle;
	\fill[left color = green!40!blue,right color=green!40!blue!60] (0,0) -- +(5,21) -- +(26,21) -- +(21,0) -- cycle;
	\foreach \x in {0,3,...,21}  \draw[-,white] (\x,0) -- +(5,21);
	\foreach \y in {0,3,...,21}  \draw[-,white] (0.25*\y,\y) -- +(21,0);	
	\draw[thick] (16,0) -- (N8);
\end{scope}

% Carro electrico
\begin{scope}[ yshift=30, xshift=110,gray, scale=0.5]
 	\draw[gray,fill=gray!30] (10,0) to[out=170, in=0] (1,1) 
 	                             to[out=180, in=30] (-3,0) -- (-10,-3)            
 	                             to[out=200,in=100] (-14,-6)
 	                             to [out=230, in=100] (-15,-11)
 	                             -- (12,-11)
 	                             to [out=60, in=-80] (11,-5) -- cycle;
    \draw[fill=white]  (-11,-3.8) to[out=200,in=100] (-14,-6) -- (-12,-5.5) to [out=120, in=210] (-11,-3.8); 	
   \draw[left color=black, right color = green!55!blue!40] (-3,0) -- (-10,-3) to (-8,-3) to (-2.5,0) -- cycle;
   \draw[left color=black, right color = green!55!blue!40] (5.5,-4) -- (-8,-4.9) -- (-8,-4) to[out=30,in=190] (0,0) to[rounded corners] (5,0) -- cycle;
   \draw[] (5.5,-4) to[out=270, in=50] (2,-11) -- (-5.7,-11) to[out=120, in=260] (-7,-5) -- cycle;
   \draw[] (-8,-6) -- + (16,1);
   \draw[gray, fill=gray!60] (3.5,-5.5) ellipse (1 and 0.5);
   \fill[black] (-7,-5) -- (-10,-5.2) -- (-6.5,-2.9) -- cycle;
   \draw[fill=gray!30] (-7,-5) -- (-7.5,-5) arc(270:90:0.9) -- (-6.5,-3.1) -- cycle;
   \fill[opacity =0.5] (-10,-11) rectangle +(18,3);
   \fill[black] (-7,-11) rectangle  +(11,0.5);
   \fill[gray!30] (-7,-10.5) rectangle  +(11,0.3);
   \draw[] (-7.5,-8) arc(30:140:4);
   \draw[] (11.5,-8) arc(30:140:4);
   \draw[left color=black, right color = green!55!blue!40] (10.1,-1) -- +(-2,0.5) -- +(-0.5,-3) -- +(0.8,-3) -- cycle;
   \draw[black, fill=red] (10.6,-4) -- +(-1.3,0) -- +(-1.0,-1) -- +(0.3,-1) arc (-90:90:0.5);
   \draw[black, fill=orange] (10.9,-5) -- +(-1.3,0) -- +(-1.0,-1) -- +(0.3,-1) arc (-90:90:0.5);
   \fill[white, opacity =0.5] (-11,-10) circle (4);
   \fill[white, opacity =0.3] (8,-10) circle (4); 	
   \fill[gray!50] (-1,-13) ellipse (14 and 0.5);
   \fill[black] (-11,-10) circle (3);
   \fill[gray!30] (-11,-10) circle (2);
   \fill[black] (-11,-10) -- +(30:2.5) -- +(0:2.5) -- cycle;
   \fill[black] (-11,-10) -- +(90:2.5) -- +(60:2.5) -- cycle;
   \fill[black] (-11,-10) -- +(150:2.5) -- +(120:2.5) -- cycle;
   \fill[black] (-11,-10) -- +(210:2.5) -- +(180:2.5) -- cycle;
   \fill[black] (-11,-10) -- +(240:2.5) -- +(270:2.5) -- cycle;
   \fill[black] (-11,-10) -- +(300:2.5) -- +(330:2.5) -- cycle;
   \fill[gray!30] (-11,-10) circle (1);
   \draw[gray!30] (-11,-10) circle (2); 	
   \fill[black] (8,-10) circle (3);
   \fill[gray!30] (8,-10) circle (2);
   \fill[black] (8,-10) -- +(30:2.5) -- +(0:2.5) -- cycle;
   \fill[black] (8,-10) -- +(90:2.5) -- +(60:2.5) -- cycle;
   \fill[black] (8,-10) -- +(150:2.5) -- +(120:2.5) -- cycle;
   \fill[black] (8,-10) -- +(210:2.5) -- +(180:2.5) -- cycle;
   \fill[black] (8,-10) -- +(240:2.5) -- +(270:2.5) -- cycle;
   \fill[black] (8,-10) -- +(300:2.5) -- +(330:2.5) -- cycle;
   \fill[gray!30] (8,-10) circle (1);
   \draw[gray!30] (8,-10) circle (2);
   \draw[-,thick, black] (2.5,-10) -- (N3); 
\end{scope}

% Bateria
\begin{scope}[xshift = 220, yshift = 15, scale=0.5]
	\fill[gray!50] (3,0) -- +(14,0) -- +(20,8)  to[out=100, in = 60] +(6,8)  -- cycle;
	\fill [black!80] (3,0) rectangle +(14,15);
	\fill [orange] (3,15) rectangle +(14,-5);
	\fill [orange] (10,15) ellipse (7 and 2);
	\fill [orange] (10,10) ellipse (7 and 2);
	\fill [gray!40] (10,15) ellipse (6 and 1.5);
	\fill [black!80] (10,0) ellipse (7 and 2);
	\fill [gray] (7,15) rectangle +(6,2);
	\draw [gray, fill=gray!50, thick] (10,17) ellipse (3 and 1);
	\fill [gray] (10,15) ellipse (3 and 1);	
	\fill[white, opacity=0.2] (10,-2) -- +(-7,0) -- +(-7,17) -- cycle;
	\draw[-,thick] (10,0) -- (N9);
\end{scope}

% turbina
\begin{scope}[xshift=190,black!80]
    \node(b) at (0,0) {};
    \draw[rotate = -15, fill=gray!30](b)+(5,0) ellipse (5 and 0.5);
    \draw[rotate = 105, fill=gray!30](b)+(5,0) ellipse (5 and 0.5);
    \draw[rotate = 225, fill=gray!30](b)+(5,0) ellipse (5 and 0.5);
    \draw[-, fill=gray!30] (b)+(-0.5,0) -- +(-0.3,0) -- +(-0.6,-22) -- +(0.6,-22) -- +(0.3,0);
    \draw[-, fill=gray!30] (b) circle (1);
    \draw[-,thick] (0,-20) -- (N5);
\end{scope}

\draw[-,thick,blue!70!green] (-5,5) rectangle +(10,10);
\draw[-,thick,blue!70!green] (-5,5) -- +(10,10);
\draw[-,thick,blue!70!green] (0,5) -- (N1);
\node[blue!80!green] at (-2,12) {AC};
\node[blue!80!green] at ( 2,7) {DC};

\node at (7,20) [text width=70] {main substation converter};
 
\end{tikzpicture}
\caption{Example of an LVDC system for urban area appliactions}
\label{fig:distribucion}
\end{figure}

An LVDC grid consists of a bidirectional AC/DC converter placed in the main substation to which it is connected different loads and generators as depicted in Fig \ref{fig:distribucion}. Different elements can be connected to an LVDC grid such as renewable energy resources, energy storage, electric vehicles and controlled loads. These elements are integrated to the grid through a power electronic converter (i.e a constant power terminal).  Consequently, the model of the LVDC grid is non-linear and requires a power flow study.  

The existence and uniqueness of the solution are, obviously, \textit{sine qua non} conditions for rigorous analysis of the stationary state of a grid and for determining an equilibrium point in small signal stability studies \cite{Review,pf_mit}. These are characteristics of the set of algebraic equations and not of the method used to find a solution. However, it is often difficult to determine if a solution of a set of non-linear equations, such as those of the power flow, is unique. A non-linear problem could give several solutions, and in some cases, the solution may not even exist. Uniqueness must not be taken for granted.

\subsection{dc power flow vs power flow in LVDC grids}

It is important to emphasize that power flow in LVDC grids is different from "dc power flow". The first is a power flow in a grid which is actually dc and incorporates constant power terminals; while the second is a linearization of the power flow equations in ac grids which, due to a pedagogic analogy, is named in this way.

\subsection{Brief state of the art}

There is an increasing interest in LVDC grids and related subjects such as dc microgrids and dc distribution. Several studies have been done about the feasibility of these technologies. For instance, \cite{overview} presented a complete description of the potentialities of LVDC grids as well as their challenges.  Potential pathways for increased use of dc technology in buildings was considered in \cite{local_dc_distribution}. A more practical approach was presented in \cite{LVDC_distribution} where a case study for a large distribution network was considered.  

Power flow analysis in LVDC grids has been presented as an extension to well known methodologies for ac grids such as Newton-Raphson or Gauss-Seidel \cite{pf_js}.  Power flow sensitivities have been also studied in \cite{sensitivities}. However, available studies in the literature are based on numerical performance but there are no theoretical studies about uniqueness of the solution. In these studies, uniqueness is taken for granted without mathematical demonstration in spite of the fact that a non-linear problem could give several solutions. To the best of the author's knowledge, this problem has not been addressed in LVDC grids \footnote{The problem has not been fully addressed in ac grids either. A result for LVDC grids could  give an insight about general ac grids}. 

\subsection{Contribution and scope}

This paper demonstrates the existence and uniqueness of the solution of the power flow in LVDC grids. This result is general since: 1) it is independent of the numerical method  2) it is independent of size and load condition of the LVDC grid and 3) it is valid for any topology of the LVDC grid. A computational simulation demonstrates the theoretical analysis using a successive approximation method. 

Comparisons of the computational performance of different algorithms is beyond of the scope of this paper in order to maintain the generality of the main result. Computational performance depends on many factors such as the implementation of the algorithm, programming language and size of the grid. 
%In spite of its apparent simplicity, the author have not been able to find this analysis in the literature in the present form.  

\subsection{Organization of the paper}

The paper is organized as follows:  Section 2 presents the basic formulation of the power flow in LVDC grids from a practical context. Next, Section 3 demonstrates the main theoretical result followed by simulations in Section 4.  Finally conclusions and references.

\section{Power Flow in LVDC grids}

The lack of reactive power and angles in LVDC grids allows some simplifications of the mathematical formulation. Nodes are classified according to the type of control, namely: constant voltage, constant power and constant resistance.  Constant voltage terminals include the main substation converter and any converter along the grid which can maintain the voltage. Other converters in the grid must be represented as constant power terminals. These include renewable energy resources, energy storage devices and controlled loads, among others.  Constant resistance terminals are linear loads as well as step nodes (i.e. nodes without generation or load). Drop controls can be considered as a linear combination of a constant power and a constant resistance terminal.

\subsection{Mathematical formulation}

Let us consider an LVDC grid as a set of nodes represented by $\mathcal{N} = \left\{1,2,...,N \right\}$, which in turns is subdivided into three nonempty and disjoint subsets $\mathcal{N}=\left\{\mathcal{V},\mathcal{R},\mathcal{P} \right\}$ according to the type of terminal,  namely:  constant voltage ($\mathcal{V}$), constant resistance ($\mathcal{R}$)  and constant power $\mathcal{P}$.  There is usually only one constant voltage terminal but the methodology can be applied to a more general case with multiple voltage-controlled terminals. Branches  are represented as a set $\mathcal{E} = \mathcal{N}\times \mathcal{N}$ with an associated constant resistance.

Nodal voltages and currents are related by the admittance matrix $G\in \mathbb{R}^{\mathcal{N}\times \mathcal{N}}$ as follows:

\begin{equation}
	\left(\begin{array}{c} 
	   I_{\mathcal{V}} \\ 		 
		 I_{\mathcal{R}} \\ 		 
		 I_{\mathcal{P}} 
	\end{array}\right)
	= 
	\left(\begin{array}{ccc} 
	   G_{\mathcal{VV}} & G_{\mathcal{VR}} & G_{\mathcal{VP}} \\ 		 
		 G_{\mathcal{RV}} & G_{\mathcal{RR}} & G_{\mathcal{RP}} \\ 		 
		 G_{\mathcal{PV}} & G_{\mathcal{PR}} & G_{\mathcal{PP}} 
	\end{array}\right)\cdot
	\left(\begin{array}{c} 
	   V_{\mathcal{V}} \\ 		 
		 V_{\mathcal{R}} \\ 		 
		 V_{\mathcal{P}} 
	\end{array}\right)
\end{equation}

In this case, $V_{\mathcal{V}}$ is known and $I_{\mathcal{R}}$ is given by (\ref{eq:d})

\begin{equation}
	I_{\mathcal{R}} = - D_{\mathcal{RR}}\cdot V_{\mathcal{R}}
	\label{eq:d}
\end{equation}

with $D_{\mathcal{RR}}$ a diagonal matrix that includes admittances of constant power terminals. Notice this matrix can be singular (e.g. in the case of step nodes).  Equation (\ref{eq:d}) is used to reduce the size of the set of algebraic equations:

\begin{equation}
	V_{\mathcal{R}} = -(D_{\mathcal{RR}}+G_{\mathcal{RR}})^{-1} \cdot (G_{\mathcal{RV}}\cdot V_{\mathcal{V}} + G_{\mathcal{RP}}\cdot V_{\mathcal{P}})
\end{equation}

Power-controlled terminals are associated with the following non-linear equation

\begin{equation}
 P_{\mathcal{P}} = diag(V_{\mathcal{P}}) \cdot I_{\mathcal{P}}
\end{equation}

Which in turn can be written as follows

\begin{equation}
 P_{\mathcal{P}} = diag(V_{\mathcal{P}}) \cdot (J_{\mathcal{P}} + B_{\mathcal{PP}} \cdot V_{\mathcal{P}})
\end{equation}

with

\begin{eqnarray*}
	J_{\mathcal{P}} = (G_{\mathcal{PV}} - G_{\mathcal{PR}}\cdot(D_{\mathcal{RR}}+G_{\mathcal{RR}})^{-1}\cdot G_{\mathcal{RV}})\cdot V_{\mathcal{V}} \\
	B_{\mathcal{PP}} = G_{\mathcal{PP}} - G_{\mathcal{PR}}\cdot(D_{\mathcal{RR}}+G_{\mathcal{RR}})^{-1}\cdot G_{\mathcal{RP}} 
\end{eqnarray*}

Therefore, the state of the LVDC grid can be completely established by solving (\ref{eq:eqnle}). 

\begin{equation}
	V_{\mathcal{P}} = B_{\mathcal{PP}}^{-1} \cdot (diag(V_{\mathcal{P}}^{-1})\cdot 
	P_{\mathcal{P}}-J_{\mathcal{P}})
	\label{eq:eqnle}
\end{equation}

In order to analyze (\ref{eq:eqnle}), let us define a map $T:\mathbb{R}^{\mathcal{P}}\rightarrow\mathbb{R}^{\mathcal{P}}$ as given in (\ref{eq:map}):

\begin{equation}
	T(V_{\mathcal{P}}) = B_{\mathcal{PP}}^{-1} \cdot (diag(V_{\mathcal{P}}^{-1})\cdot 
	P_{\mathcal{P}}-J_{\mathcal{P}})
	\label{eq:map}
\end{equation}

Notice that T is a non-linear map and hence, uniqueness of the solution must not be taken for granted.

\subsection{Practical considerations}

Let us consider the following few practical assumptions
\begin{description}
 \item[A1] the graph is connected (i.e. there are no islands in the feeder). 
 \item[A2] there is at least one constant power terminal and one constant voltage terminal
 \item[A3] feasible voltages remain in a given interval 
\begin{equation}
0<v_{min}\leq V\leq v_{max}
\end{equation}
 \item[A4] short circuit currents are higher than normal operation currents for all constant power terminals.
\end{description}

Each of these assumptions is completely justified in real power system applications.
(A1-A2) guarantee $B_{\mathcal{PP}}$ is non-singular. (A3) is required for voltage regulation and for physical constraints in the converters.  Finally (A4) is an obvious yet useful observation for any electric system.

\section{Convergence Analysis}

We now analyze (\ref{eq:map}) in order to determine the existence and uniqueness of the solution. To do this, we must demonstrate that $T$ is a contraction mapping defined as follows\footnote{This definition can be extended to general Banch spaces.  However, a version in $\mathbb{R}^n$ is enough for the purposes of this paper.}:

\begin{definition}[contraction mapping]\ 

Let $\mathcal{B}=\left\{x: \left\|x\right\|\leq r \right\}$ be a closed ball in $\mathbb{R}^n$, and let $T:\mathcal{B}\rightarrow \mathbb{R}^n$. Then $T$ is said to be a contraction mapping if there is an $\alpha$ such that $\left\|T(x)-T(y) \right\|\leq \alpha\left\|x-y \right\|$, with $0\leq \alpha < 1, \; \forall\; x,y \in \mathcal{B}$
\end{definition}

Now we can present our main result:

\begin{proposition}[Contraction of the power flow]\ \\
  An LVDC grid represented by (\ref{eq:eqnle}) with the assumptions (A1-A4) has a unique solution which can be obtained by the method of successive approximations with the map (\ref{eq:map}) and contraction constant given by (\ref{eq:condicion})
	\begin{equation}
		\alpha = \frac{\left\| B_{\mathcal{PP}}^{-1}\right\|\cdot \left\| P_{\mathcal{P}}\right\|}{v_{min}^2}
		\label{eq:condicion}
	\end{equation}
	
\end{proposition}

\textbf{proff:} In order to prof this proposition, we use the contraction mapping theorem (see \cite{sholomo} for details) which states that if $T$ is a contraction mapping in $\mathcal{B}$ then there is a \textit{unique} vector $x_0 \in \mathcal{B}$ satisfying $x_0=T(x_0)$.

Select two different values of voltages $V_{\mathcal{P}}$ and $U_{\mathcal{P}}$ in the constant power terminals as follows:
\begin{eqnarray*}
	\left\| T(V_{\mathcal{P}})-T(U_{\mathcal{P}}) \right\|  = 	\\
	\left\| B_{\mathcal{PP}}^{-1}\cdot diag(1/V_{\mathcal{P}}-1/U_{\mathcal{P}}) \cdot P_{\mathcal{P}} \right\|	\\
	\leq \left\| B_{\mathcal{PP}}^{-1}\right\|\cdot\left\|diag(1/V_{\mathcal{P}}-1/U_{\mathcal{P}})\right\|\cdot \left\|P_{\mathcal{P}} \right\| \\
	\leq \left(\frac{\left\| B_{\mathcal{PP}}^{-1}\right\|\cdot \left\| P_{\mathcal{P}}\right\|}{v_{min}^2}	\right)\cdot \left\| V_{\mathcal{P}}-U_{\mathcal{P}}\right\| \\
	= \alpha \cdot \left\| V_{\mathcal{P}}-U_{\mathcal{P}}\right\|
\end{eqnarray*}

where $v_{min}$ is the minimum voltage according to assumption (A3). It only remains to establish if the constant $\alpha$ is lower than 1. Let us use a matrix norm as

\begin{equation}
  \left\| X \right\| = \underset{ij}{\text{max}}\;\left\{ |x_{ij}| \right\}
\end{equation}

since $B_{\mathcal{PP}}^{-1}$ is diagonal dominant, then 

\begin{equation}
\left\| B_{\mathcal{PP}}^{-1}\right\| = max\left\{ r_{kk} \right\}
\end{equation}

where $r_{kk}$ is the Thevening impedance in each node (i.e the element $k$ on the diagonal of the matrix $B_{\mathcal{PP}}^{-1}$).  Then $\alpha$ can be expressed as

\begin{equation}
	\alpha = \underset{k\in \mathcal{P}}{\text{max}}\left\{\frac{P_{k}/v_{min}}{v_{min}/r_{kk}}\right\} 
	\label{eq:cconstant}
\end{equation}

which is the ratio between operational and short circuit currents at minimum voltage.  This value is lower than one due to (A4); hence, the prof is completed.

\textbf{Remark 3.1:} Notice that (\ref{eq:condicion}) can be directly evaluated before the power flow calculation.  This condition is inherent in the system and not in the computational implementation.

\textbf{Remark 3.2:} The theorem guarantees uniqueness of the solution in $\mathcal{B}$. As aforementioned, it does not depend on the implementation; any algorithm that achieves convergence will find a point in $\mathcal{B}$.

\subsection{Successive approximations}

Several methodologies from ac systems can be adapted to LVDC grids.  Many of these are based on the classic Gauss-Seidel and/or Newton-Raphson methods. No method is guaranteed to be faster than the other in every case. Here, a successive approximation is used since it can be directly obtained from the map $T$; this method applies iteratively the map $T$ until achieving convergence as follows:

\begin{equation}
	V_{\mathcal{P}(k+1)} = T(V_{\mathcal{P}(k)})
\end{equation}

where the sub-index $k$ represents the iteration.  This methodology is classic in the power systems literature. In fact, the Gauss-Seidel method is just a small modification of this principle (in this method the values of $V_{\mathcal{P}(k)}$ obtained in the k-th iteration remain unchanged during the entire iteration while the Gauss-Seidel method uses the new values as soon as they are obtained). In addition, the backward-forward sweep algorithm can be interpreted as a computationally efficient implementation of the same principle. Therefore, the analysis of these algorithms is basically equivalent.

\section{Computational results}

A power flow was evaluated in the medium voltage dc distribution system shown in Fig \ref{fig:distribucion} with parameters given in Table \ref{tab:parametros}.  Convergence of the method was analyzed under different initial conditions for $\mathcal{B}=\left\{V \in \mathbb{R}^{\mathcal{P}}: 0.55<V_{i}<1.5\right\}$.  The contraction constant was calculated using (\ref{eq:cconstant}) as $\alpha=0.00475$.  The algorithm converged in less than 5 iterations  regardless  the initial point.

\begin{table}[bt]
\centering
\caption{Parameters of the LVDC}
\label{tab:parametros}
\begin{tabular}{|c|c|c|c|c|}
\hline
From & To & $r(pu)$ & Type & $P/R(pu)$\\
\hline\hline
    1  &  2  &  0.0050 & step-node &      \\
    2  &  3  &  0.0015 & P 				& -0.8 \\
    2  &  4  &  0.0020 & P					& -1.3 \\
    4  &  5  &  0.0018 & P					&  0.5 \\
    2  &  6  &  0.0023 & R 				&  2.0 \\
    6  &  7  &  0.0017 & step-node &  \\
    7  &  8  &  0.0021 & P					&  0.3 \\
    7  &  9  &  0.0013 & P					& -0.7 \\
    3  & 10  &  0.0015 & R 				&  1.25 \\
\hline	
\end{tabular}
\end{table}

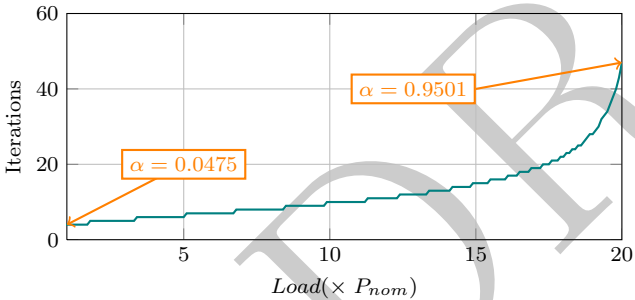
\begin{figure}[htb]
\footnotesize
\centering

\begin{tikzpicture}
\begin{axis}[scale only axis,width=7.3cm, height=3.0cm,xmajorgrids,ymajorgrids, xlabel={$Load (\times \; P_{nom})$}, ylabel ={Iterations}, ymax = 60, ymin=0, xmin=1, xmax = 20]
\addplot [thick, blue!50!green]
coordinates{( 0.100000 , 3.000000 )
( 0.200000 , 3.000000 )
( 0.300000 , 3.000000 )
( 0.400000 , 3.000000 )
( 0.500000 , 3.000000 )
( 0.600000 , 3.000000 )
( 0.700000 , 4.000000 )
( 0.800000 , 4.000000 )
( 0.900000 , 4.000000 )
( 1.000000 , 4.000000 )
( 1.100000 , 4.000000 )
( 1.200000 , 4.000000 )
( 1.300000 , 4.000000 )
( 1.400000 , 4.000000 )
( 1.500000 , 4.000000 )
( 1.600000 , 4.000000 )
( 1.700000 , 4.000000 )
( 1.800000 , 5.000000 )
( 1.900000 , 5.000000 )
( 2.000000 , 5.000000 )
( 2.100000 , 5.000000 )
( 2.200000 , 5.000000 )
( 2.300000 , 5.000000 )
( 2.400000 , 5.000000 )
( 2.500000 , 5.000000 )
( 2.600000 , 5.000000 )
( 2.700000 , 5.000000 )
( 2.800000 , 5.000000 )
( 2.900000 , 5.000000 )
( 3.000000 , 5.000000 )
( 3.100000 , 5.000000 )
( 3.200000 , 5.000000 )
( 3.300000 , 5.000000 )
( 3.400000 , 6.000000 )
( 3.500000 , 6.000000 )
( 3.600000 , 6.000000 )
( 3.700000 , 6.000000 )
( 3.800000 , 6.000000 )
( 3.900000 , 6.000000 )
( 4.000000 , 6.000000 )
( 4.100000 , 6.000000 )
( 4.200000 , 6.000000 )
( 4.300000 , 6.000000 )
( 4.400000 , 6.000000 )
( 4.500000 , 6.000000 )
( 4.600000 , 6.000000 )
( 4.700000 , 6.000000 )
( 4.800000 , 6.000000 )
( 4.900000 , 6.000000 )
( 5.000000 , 6.000000 )
( 5.100000 , 7.000000 )
( 5.200000 , 7.000000 )
( 5.300000 , 7.000000 )
( 5.400000 , 7.000000 )
( 5.500000 , 7.000000 )
( 5.600000 , 7.000000 )
( 5.700000 , 7.000000 )
( 5.800000 , 7.000000 )
( 5.900000 , 7.000000 )
( 6.000000 , 7.000000 )
( 6.100000 , 7.000000 )
( 6.200000 , 7.000000 )
( 6.300000 , 7.000000 )
( 6.400000 , 7.000000 )
( 6.500000 , 7.000000 )
( 6.600000 , 7.000000 )
( 6.700000 , 7.000000 )
( 6.800000 , 8.000000 )
( 6.900000 , 8.000000 )
( 7.000000 , 8.000000 )
( 7.100000 , 8.000000 )
( 7.200000 , 8.000000 )
( 7.300000 , 8.000000 )
( 7.400000 , 8.000000 )
( 7.500000 , 8.000000 )
( 7.600000 , 8.000000 )
( 7.700000 , 8.000000 )
( 7.800000 , 8.000000 )
( 7.900000 , 8.000000 )
( 8.000000 , 8.000000 )
( 8.100000 , 8.000000 )
( 8.200000 , 8.000000 )
( 8.300000 , 8.000000 )
( 8.400000 , 8.000000 )
( 8.500000 , 9.000000 )
( 8.600000 , 9.000000 )
( 8.700000 , 9.000000 )
( 8.800000 , 9.000000 )
( 8.900000 , 9.000000 )
( 9.000000 , 9.000000 )
( 9.100000 , 9.000000 )
( 9.200000 , 9.000000 )
( 9.300000 , 9.000000 )
( 9.400000 , 9.000000 )
( 9.500000 , 9.000000 )
( 9.600000 , 9.000000 )
( 9.700000 , 9.000000 )
( 9.800000 , 9.000000 )
( 9.900000 , 10.000000 )
( 10.000000 , 10.000000 )
( 10.100000 , 10.000000 )
( 10.200000 , 10.000000 )
( 10.300000 , 10.000000 )
( 10.400000 , 10.000000 )
( 10.500000 , 10.000000 )
( 10.600000 , 10.000000 )
( 10.700000 , 10.000000 )
( 10.800000 , 10.000000 )
( 10.900000 , 10.000000 )
( 11.000000 , 10.000000 )
( 11.100000 , 10.000000 )
( 11.200000 , 10.000000 )
( 11.300000 , 11.000000 )
( 11.400000 , 11.000000 )
( 11.500000 , 11.000000 )
( 11.600000 , 11.000000 )
( 11.700000 , 11.000000 )
( 11.800000 , 11.000000 )
( 11.900000 , 11.000000 )
( 12.000000 , 11.000000 )
( 12.100000 , 11.000000 )
( 12.200000 , 11.000000 )
( 12.300000 , 11.000000 )
( 12.400000 , 12.000000 )
( 12.500000 , 12.000000 )
( 12.600000 , 12.000000 )
( 12.700000 , 12.000000 )
( 12.800000 , 12.000000 )
( 12.900000 , 12.000000 )
( 13.000000 , 12.000000 )
( 13.100000 , 12.000000 )
( 13.200000 , 12.000000 )
( 13.300000 , 12.000000 )
( 13.400000 , 13.000000 )
( 13.500000 , 13.000000 )
( 13.600000 , 13.000000 )
( 13.700000 , 13.000000 )
( 13.800000 , 13.000000 )
( 13.900000 , 13.000000 )
( 14.000000 , 13.000000 )
( 14.100000 , 13.000000 )
( 14.200000 , 14.000000 )
( 14.300000 , 14.000000 )
( 14.400000 , 14.000000 )
( 14.500000 , 14.000000 )
( 14.600000 , 14.000000 )
( 14.700000 , 14.000000 )
( 14.800000 , 14.000000 )
( 14.900000 , 15.000000 )
( 15.000000 , 15.000000 )
( 15.100000 , 15.000000 )
( 15.200000 , 15.000000 )
( 15.300000 , 15.000000 )
( 15.400000 , 15.000000 )
( 15.500000 , 16.000000 )
( 15.600000 , 16.000000 )
( 15.700000 , 16.000000 )
( 15.800000 , 16.000000 )
( 15.900000 , 16.000000 )
( 16.000000 , 16.000000 )
( 16.100000 , 17.000000 )
( 16.200000 , 17.000000 )
( 16.300000 , 17.000000 )
( 16.400000 , 17.000000 )
( 16.500000 , 18.000000 )
( 16.600000 , 18.000000 )
( 16.700000 , 18.000000 )
( 16.800000 , 18.000000 )
( 16.900000 , 19.000000 )
( 17.000000 , 19.000000 )
( 17.100000 , 19.000000 )
( 17.200000 , 19.000000 )
( 17.300000 , 20.000000 )
( 17.400000 , 20.000000 )
( 17.500000 , 20.000000 )
( 17.600000 , 21.000000 )
( 17.700000 , 21.000000 )
( 17.800000 , 21.000000 )
( 17.900000 , 22.000000 )
( 18.000000 , 22.000000 )
( 18.100000 , 23.000000 )
( 18.200000 , 23.000000 )
( 18.300000 , 24.000000 )
( 18.400000 , 24.000000 )
( 18.500000 , 25.000000 )
( 18.600000 , 25.000000 )
( 18.700000 , 26.000000 )
( 18.800000 , 27.000000 )
( 18.900000 , 28.000000 )
( 19.000000 , 28.000000 )
( 19.100000 , 29.000000 )
( 19.200000 , 30.000000 )
( 19.300000 , 32.000000 )
( 19.400000 , 33.000000 )
( 19.500000 , 34.000000 )
( 19.600000 , 36.000000 )
( 19.700000 , 38.000000 )
( 19.800000 , 40.000000 )
( 19.900000 , 43.000000 )
( 20.000000 , 47.000000 )};
\draw[<-,thick,orange] (axis cs:1,4) -- (axis cs:5,20) node [draw,fill=white] {$\alpha=0.0475$};
\draw[<-,thick,orange] (axis cs:20,47) -- (axis cs:15,40) node [draw,fill=white,left] {$\alpha=0.9501$};
\end{axis}
\end{tikzpicture}
\caption{Convergency properties for different load conditions}
\label{fig:cargabilidad}
\end{figure}

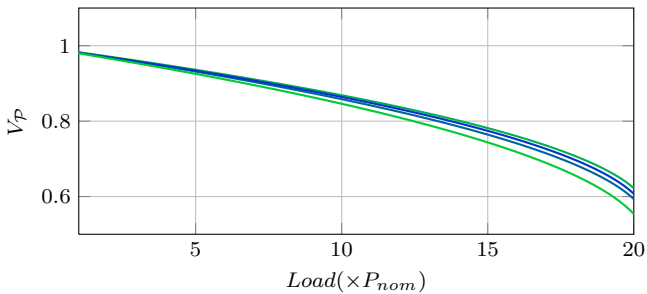
\begin{figure}[tb]
\footnotesize
\centering
\begin{tikzpicture}
\begin{axis}[scale only axis,width=7.3cm, height=3.0cm,xmajorgrids,ymajorgrids, xlabel={$Load (\times P_{nom})$}, ylabel ={$V_{\mathcal{P}}$}, ymax = 1.1, ymin=0.5, xmin=1, xmax = 20]
\addplot [thick, blue!60!green]coordinates{( 0.100000 , 0.992392 )
( 0.200000 , 0.991227 )
( 0.300000 , 0.990058 )
( 0.400000 , 0.988887 )
( 0.500000 , 0.987713 )
( 0.600000 , 0.986536 )
( 0.700000 , 0.985356 )
( 0.800000 , 0.984174 )
( 0.900000 , 0.982988 )
( 1.000000 , 0.981799 )
( 1.100000 , 0.980607 )
( 1.200000 , 0.979412 )
( 1.300000 , 0.978214 )
( 1.400000 , 0.977013 )
( 1.500000 , 0.975808 )
( 1.600000 , 0.974601 )
( 1.700000 , 0.973390 )
( 1.800000 , 0.972176 )
( 1.900000 , 0.970959 )
( 2.000000 , 0.969738 )
( 2.100000 , 0.968515 )
( 2.200000 , 0.967288 )
( 2.300000 , 0.966057 )
( 2.400000 , 0.964823 )
( 2.500000 , 0.963586 )
( 2.600000 , 0.962345 )
( 2.700000 , 0.961101 )
( 2.800000 , 0.959854 )
( 2.900000 , 0.958602 )
( 3.000000 , 0.957348 )
( 3.100000 , 0.956089 )
( 3.200000 , 0.954828 )
( 3.300000 , 0.953562 )
( 3.400000 , 0.952293 )
( 3.500000 , 0.951020 )
( 3.600000 , 0.949743 )
( 3.700000 , 0.948463 )
( 3.800000 , 0.947178 )
( 3.900000 , 0.945890 )
( 4.000000 , 0.944598 )
( 4.100000 , 0.943302 )
( 4.200000 , 0.942002 )
( 4.300000 , 0.940699 )
( 4.400000 , 0.939391 )
( 4.500000 , 0.938079 )
( 4.600000 , 0.936763 )
( 4.700000 , 0.935443 )
( 4.800000 , 0.934118 )
( 4.900000 , 0.932790 )
( 5.000000 , 0.931457 )
( 5.100000 , 0.930120 )
( 5.200000 , 0.928778 )
( 5.300000 , 0.927433 )
( 5.400000 , 0.926082 )
( 5.500000 , 0.924728 )
( 5.600000 , 0.923368 )
( 5.700000 , 0.922005 )
( 5.800000 , 0.920636 )
( 5.900000 , 0.919263 )
( 6.000000 , 0.917885 )
( 6.100000 , 0.916503 )
( 6.200000 , 0.915116 )
( 6.300000 , 0.913723 )
( 6.400000 , 0.912326 )
( 6.500000 , 0.910924 )
( 6.600000 , 0.909517 )
( 6.700000 , 0.908105 )
( 6.800000 , 0.906688 )
( 6.900000 , 0.905266 )
( 7.000000 , 0.903838 )
( 7.100000 , 0.902406 )
( 7.200000 , 0.900967 )
( 7.300000 , 0.899524 )
( 7.400000 , 0.898075 )
( 7.500000 , 0.896620 )
( 7.600000 , 0.895160 )
( 7.700000 , 0.893694 )
( 7.800000 , 0.892223 )
( 7.900000 , 0.890745 )
( 8.000000 , 0.889262 )
( 8.100000 , 0.887773 )
( 8.200000 , 0.886278 )
( 8.300000 , 0.884777 )
( 8.400000 , 0.883269 )
( 8.500000 , 0.881756 )
( 8.600000 , 0.880236 )
( 8.700000 , 0.878710 )
( 8.800000 , 0.877177 )
( 8.900000 , 0.875638 )
( 9.000000 , 0.874092 )
( 9.100000 , 0.872540 )
( 9.200000 , 0.870980 )
( 9.300000 , 0.869414 )
( 9.400000 , 0.867841 )
( 9.500000 , 0.866261 )
( 9.600000 , 0.864673 )
( 9.700000 , 0.863078 )
( 9.800000 , 0.861476 )
( 9.900000 , 0.859867 )
( 10.000000 , 0.858249 )
( 10.100000 , 0.856625 )
( 10.200000 , 0.854992 )
( 10.300000 , 0.853351 )
( 10.400000 , 0.851703 )
( 10.500000 , 0.850046 )
( 10.600000 , 0.848381 )
( 10.700000 , 0.846708 )
( 10.800000 , 0.845026 )
( 10.900000 , 0.843335 )
( 11.000000 , 0.841636 )
( 11.100000 , 0.839927 )
( 11.200000 , 0.838210 )
( 11.300000 , 0.836483 )
( 11.400000 , 0.834747 )
( 11.500000 , 0.833002 )
( 11.600000 , 0.831247 )
( 11.700000 , 0.829482 )
( 11.800000 , 0.827707 )
( 11.900000 , 0.825921 )
( 12.000000 , 0.824126 )
( 12.100000 , 0.822320 )
( 12.200000 , 0.820503 )
( 12.300000 , 0.818675 )
( 12.400000 , 0.816836 )
( 12.500000 , 0.814986 )
( 12.600000 , 0.813124 )
( 12.700000 , 0.811251 )
( 12.800000 , 0.809365 )
( 12.900000 , 0.807468 )
( 13.000000 , 0.805558 )
( 13.100000 , 0.803635 )
( 13.200000 , 0.801699 )
( 13.300000 , 0.799750 )
( 13.400000 , 0.797787 )
( 13.500000 , 0.795811 )
( 13.600000 , 0.793820 )
( 13.700000 , 0.791816 )
( 13.800000 , 0.789796 )
( 13.900000 , 0.787761 )
( 14.000000 , 0.785712 )
( 14.100000 , 0.783646 )
( 14.200000 , 0.781564 )
( 14.300000 , 0.779466 )
( 14.400000 , 0.777351 )
( 14.500000 , 0.775218 )
( 14.600000 , 0.773068 )
( 14.700000 , 0.770900 )
( 14.800000 , 0.768713 )
( 14.900000 , 0.766507 )
( 15.000000 , 0.764281 )
( 15.100000 , 0.762035 )
( 15.200000 , 0.759768 )
( 15.300000 , 0.757480 )
( 15.400000 , 0.755170 )
( 15.500000 , 0.752838 )
( 15.600000 , 0.750482 )
( 15.700000 , 0.748102 )
( 15.800000 , 0.745697 )
( 15.900000 , 0.743267 )
( 16.000000 , 0.740811 )
( 16.100000 , 0.738327 )
( 16.200000 , 0.735815 )
( 16.300000 , 0.733274 )
( 16.400000 , 0.730703 )
( 16.500000 , 0.728101 )
( 16.600000 , 0.725466 )
( 16.700000 , 0.722797 )
( 16.800000 , 0.720094 )
( 16.900000 , 0.717354 )
( 17.000000 , 0.714576 )
( 17.100000 , 0.711758 )
( 17.200000 , 0.708899 )
( 17.300000 , 0.705997 )
( 17.400000 , 0.703049 )
( 17.500000 , 0.700054 )
( 17.600000 , 0.697009 )
( 17.700000 , 0.693911 )
( 17.800000 , 0.690759 )
( 17.900000 , 0.687548 )
( 18.000000 , 0.684275 )
( 18.100000 , 0.680937 )
( 18.200000 , 0.677529 )
( 18.300000 , 0.674047 )
( 18.400000 , 0.670486 )
( 18.500000 , 0.666840 )
( 18.600000 , 0.663102 )
( 18.700000 , 0.659266 )
( 18.800000 , 0.655323 )
( 18.900000 , 0.651264 )
( 19.000000 , 0.647077 )
( 19.100000 , 0.642749 )
( 19.200000 , 0.638266 )
( 19.300000 , 0.633609 )
( 19.400000 , 0.628756 )
( 19.500000 , 0.623680 )
( 19.600000 , 0.618345 )
( 19.700000 , 0.612709 )
( 19.800000 , 0.606712 )
( 19.900000 , 0.600275 )
( 20.000000 , 0.593281 )
};
\addplot [thick, blue!30!green]coordinates{( 0.100000 , 0.992482 )
( 0.200000 , 0.991408 )
( 0.300000 , 0.990331 )
( 0.400000 , 0.989251 )
( 0.500000 , 0.988169 )
( 0.600000 , 0.987083 )
( 0.700000 , 0.985995 )
( 0.800000 , 0.984905 )
( 0.900000 , 0.983811 )
( 1.000000 , 0.982715 )
( 1.100000 , 0.981615 )
( 1.200000 , 0.980513 )
( 1.300000 , 0.979408 )
( 1.400000 , 0.978301 )
( 1.500000 , 0.977190 )
( 1.600000 , 0.976076 )
( 1.700000 , 0.974959 )
( 1.800000 , 0.973840 )
( 1.900000 , 0.972717 )
( 2.000000 , 0.971591 )
( 2.100000 , 0.970462 )
( 2.200000 , 0.969330 )
( 2.300000 , 0.968195 )
( 2.400000 , 0.967057 )
( 2.500000 , 0.965915 )
( 2.600000 , 0.964771 )
( 2.700000 , 0.963623 )
( 2.800000 , 0.962472 )
( 2.900000 , 0.961317 )
( 3.000000 , 0.960160 )
( 3.100000 , 0.958999 )
( 3.200000 , 0.957834 )
( 3.300000 , 0.956666 )
( 3.400000 , 0.955495 )
( 3.500000 , 0.954321 )
( 3.600000 , 0.953142 )
( 3.700000 , 0.951961 )
( 3.800000 , 0.950775 )
( 3.900000 , 0.949587 )
( 4.000000 , 0.948394 )
( 4.100000 , 0.947198 )
( 4.200000 , 0.945998 )
( 4.300000 , 0.944795 )
( 4.400000 , 0.943588 )
( 4.500000 , 0.942376 )
( 4.600000 , 0.941162 )
( 4.700000 , 0.939943 )
( 4.800000 , 0.938720 )
( 4.900000 , 0.937494 )
( 5.000000 , 0.936263 )
( 5.100000 , 0.935029 )
( 5.200000 , 0.933790 )
( 5.300000 , 0.932548 )
( 5.400000 , 0.931301 )
( 5.500000 , 0.930050 )
( 5.600000 , 0.928795 )
( 5.700000 , 0.927535 )
( 5.800000 , 0.926272 )
( 5.900000 , 0.925004 )
( 6.000000 , 0.923731 )
( 6.100000 , 0.922454 )
( 6.200000 , 0.921173 )
( 6.300000 , 0.919887 )
( 6.400000 , 0.918597 )
( 6.500000 , 0.917302 )
( 6.600000 , 0.916002 )
( 6.700000 , 0.914698 )
( 6.800000 , 0.913389 )
( 6.900000 , 0.912075 )
( 7.000000 , 0.910756 )
( 7.100000 , 0.909432 )
( 7.200000 , 0.908103 )
( 7.300000 , 0.906769 )
( 7.400000 , 0.905430 )
( 7.500000 , 0.904086 )
( 7.600000 , 0.902737 )
( 7.700000 , 0.901382 )
( 7.800000 , 0.900022 )
( 7.900000 , 0.898657 )
( 8.000000 , 0.897286 )
( 8.100000 , 0.895910 )
( 8.200000 , 0.894528 )
( 8.300000 , 0.893141 )
( 8.400000 , 0.891747 )
( 8.500000 , 0.890348 )
( 8.600000 , 0.888943 )
( 8.700000 , 0.887532 )
( 8.800000 , 0.886115 )
( 8.900000 , 0.884692 )
( 9.000000 , 0.883263 )
( 9.100000 , 0.881827 )
( 9.200000 , 0.880385 )
( 9.300000 , 0.878937 )
( 9.400000 , 0.877482 )
( 9.500000 , 0.876021 )
( 9.600000 , 0.874552 )
( 9.700000 , 0.873077 )
( 9.800000 , 0.871596 )
( 9.900000 , 0.870107 )
( 10.000000 , 0.868611 )
( 10.100000 , 0.867108 )
( 10.200000 , 0.865597 )
( 10.300000 , 0.864079 )
( 10.400000 , 0.862554 )
( 10.500000 , 0.861021 )
( 10.600000 , 0.859481 )
( 10.700000 , 0.857932 )
( 10.800000 , 0.856376 )
( 10.900000 , 0.854811 )
( 11.000000 , 0.853238 )
( 11.100000 , 0.851657 )
( 11.200000 , 0.850068 )
( 11.300000 , 0.848470 )
( 11.400000 , 0.846863 )
( 11.500000 , 0.845247 )
( 11.600000 , 0.843622 )
( 11.700000 , 0.841988 )
( 11.800000 , 0.840344 )
( 11.900000 , 0.838691 )
( 12.000000 , 0.837029 )
( 12.100000 , 0.835356 )
( 12.200000 , 0.833674 )
( 12.300000 , 0.831981 )
( 12.400000 , 0.830278 )
( 12.500000 , 0.828564 )
( 12.600000 , 0.826839 )
( 12.700000 , 0.825104 )
( 12.800000 , 0.823357 )
( 12.900000 , 0.821599 )
( 13.000000 , 0.819829 )
( 13.100000 , 0.818047 )
( 13.200000 , 0.816253 )
( 13.300000 , 0.814447 )
( 13.400000 , 0.812628 )
( 13.500000 , 0.810796 )
( 13.600000 , 0.808951 )
( 13.700000 , 0.807093 )
( 13.800000 , 0.805220 )
( 13.900000 , 0.803334 )
( 14.000000 , 0.801433 )
( 14.100000 , 0.799518 )
( 14.200000 , 0.797587 )
( 14.300000 , 0.795641 )
( 14.400000 , 0.793680 )
( 14.500000 , 0.791702 )
( 14.600000 , 0.789707 )
( 14.700000 , 0.787696 )
( 14.800000 , 0.785666 )
( 14.900000 , 0.783620 )
( 15.000000 , 0.781554 )
( 15.100000 , 0.779470 )
( 15.200000 , 0.777366 )
( 15.300000 , 0.775242 )
( 15.400000 , 0.773098 )
( 15.500000 , 0.770933 )
( 15.600000 , 0.768745 )
( 15.700000 , 0.766536 )
( 15.800000 , 0.764303 )
( 15.900000 , 0.762046 )
( 16.000000 , 0.759764 )
( 16.100000 , 0.757457 )
( 16.200000 , 0.755124 )
( 16.300000 , 0.752763 )
( 16.400000 , 0.750374 )
( 16.500000 , 0.747955 )
( 16.600000 , 0.745506 )
( 16.700000 , 0.743026 )
( 16.800000 , 0.740512 )
( 16.900000 , 0.737964 )
( 17.000000 , 0.735381 )
( 17.100000 , 0.732761 )
( 17.200000 , 0.730101 )
( 17.300000 , 0.727401 )
( 17.400000 , 0.724659 )
( 17.500000 , 0.721872 )
( 17.600000 , 0.719038 )
( 17.700000 , 0.716155 )
( 17.800000 , 0.713220 )
( 17.900000 , 0.710230 )
( 18.000000 , 0.707183 )
( 18.100000 , 0.704073 )
( 18.200000 , 0.700899 )
( 18.300000 , 0.697655 )
( 18.400000 , 0.694336 )
( 18.500000 , 0.690937 )
( 18.600000 , 0.687453 )
( 18.700000 , 0.683876 )
( 18.800000 , 0.680198 )
( 18.900000 , 0.676411 )
( 19.000000 , 0.672504 )
( 19.100000 , 0.668465 )
( 19.200000 , 0.664279 )
( 19.300000 , 0.659930 )
( 19.400000 , 0.655397 )
( 19.500000 , 0.650653 )
( 19.600000 , 0.645666 )
( 19.700000 , 0.640395 )
( 19.800000 , 0.634785 )
( 19.900000 , 0.628760 )
( 20.000000 , 0.622210 )
};
\addplot [thick, blue!80!green]coordinates{( 0.100000 , 0.991315 )
( 0.200000 , 0.990215 )
( 0.300000 , 0.989111 )
( 0.400000 , 0.988005 )
( 0.500000 , 0.986896 )
( 0.600000 , 0.985784 )
( 0.700000 , 0.984669 )
( 0.800000 , 0.983552 )
( 0.900000 , 0.982431 )
( 1.000000 , 0.981308 )
( 1.100000 , 0.980182 )
( 1.200000 , 0.979052 )
( 1.300000 , 0.977920 )
( 1.400000 , 0.976784 )
( 1.500000 , 0.975646 )
( 1.600000 , 0.974505 )
( 1.700000 , 0.973360 )
( 1.800000 , 0.972212 )
( 1.900000 , 0.971061 )
( 2.000000 , 0.969907 )
( 2.100000 , 0.968750 )
( 2.200000 , 0.967590 )
( 2.300000 , 0.966426 )
( 2.400000 , 0.965259 )
( 2.500000 , 0.964089 )
( 2.600000 , 0.962915 )
( 2.700000 , 0.961738 )
( 2.800000 , 0.960558 )
( 2.900000 , 0.959374 )
( 3.000000 , 0.958187 )
( 3.100000 , 0.956996 )
( 3.200000 , 0.955802 )
( 3.300000 , 0.954605 )
( 3.400000 , 0.953403 )
( 3.500000 , 0.952198 )
( 3.600000 , 0.950990 )
( 3.700000 , 0.949778 )
( 3.800000 , 0.948562 )
( 3.900000 , 0.947342 )
( 4.000000 , 0.946119 )
( 4.100000 , 0.944892 )
( 4.200000 , 0.943661 )
( 4.300000 , 0.942426 )
( 4.400000 , 0.941187 )
( 4.500000 , 0.939944 )
( 4.600000 , 0.938698 )
( 4.700000 , 0.937447 )
( 4.800000 , 0.936192 )
( 4.900000 , 0.934933 )
( 5.000000 , 0.933670 )
( 5.100000 , 0.932403 )
( 5.200000 , 0.931132 )
( 5.300000 , 0.929856 )
( 5.400000 , 0.928576 )
( 5.500000 , 0.927292 )
( 5.600000 , 0.926003 )
( 5.700000 , 0.924710 )
( 5.800000 , 0.923412 )
( 5.900000 , 0.922110 )
( 6.000000 , 0.920803 )
( 6.100000 , 0.919492 )
( 6.200000 , 0.918176 )
( 6.300000 , 0.916855 )
( 6.400000 , 0.915530 )
( 6.500000 , 0.914199 )
( 6.600000 , 0.912864 )
( 6.700000 , 0.911524 )
( 6.800000 , 0.910179 )
( 6.900000 , 0.908829 )
( 7.000000 , 0.907473 )
( 7.100000 , 0.906113 )
( 7.200000 , 0.904748 )
( 7.300000 , 0.903377 )
( 7.400000 , 0.902000 )
( 7.500000 , 0.900619 )
( 7.600000 , 0.899232 )
( 7.700000 , 0.897839 )
( 7.800000 , 0.896441 )
( 7.900000 , 0.895037 )
( 8.000000 , 0.893628 )
( 8.100000 , 0.892213 )
( 8.200000 , 0.890792 )
( 8.300000 , 0.889365 )
( 8.400000 , 0.887931 )
( 8.500000 , 0.886492 )
( 8.600000 , 0.885047 )
( 8.700000 , 0.883596 )
( 8.800000 , 0.882138 )
( 8.900000 , 0.880674 )
( 9.000000 , 0.879203 )
( 9.100000 , 0.877726 )
( 9.200000 , 0.876242 )
( 9.300000 , 0.874751 )
( 9.400000 , 0.873254 )
( 9.500000 , 0.871749 )
( 9.600000 , 0.870238 )
( 9.700000 , 0.868719 )
( 9.800000 , 0.867193 )
( 9.900000 , 0.865660 )
( 10.000000 , 0.864120 )
( 10.100000 , 0.862572 )
( 10.200000 , 0.861016 )
( 10.300000 , 0.859453 )
( 10.400000 , 0.857881 )
( 10.500000 , 0.856302 )
( 10.600000 , 0.854715 )
( 10.700000 , 0.853119 )
( 10.800000 , 0.851515 )
( 10.900000 , 0.849903 )
( 11.000000 , 0.848282 )
( 11.100000 , 0.846652 )
( 11.200000 , 0.845013 )
( 11.300000 , 0.843365 )
( 11.400000 , 0.841708 )
( 11.500000 , 0.840042 )
( 11.600000 , 0.838366 )
( 11.700000 , 0.836680 )
( 11.800000 , 0.834985 )
( 11.900000 , 0.833280 )
( 12.000000 , 0.831564 )
( 12.100000 , 0.829838 )
( 12.200000 , 0.828102 )
( 12.300000 , 0.826354 )
( 12.400000 , 0.824596 )
( 12.500000 , 0.822827 )
( 12.600000 , 0.821046 )
( 12.700000 , 0.819254 )
( 12.800000 , 0.817450 )
( 12.900000 , 0.815634 )
( 13.000000 , 0.813805 )
( 13.100000 , 0.811964 )
( 13.200000 , 0.810111 )
( 13.300000 , 0.808244 )
( 13.400000 , 0.806364 )
( 13.500000 , 0.804470 )
( 13.600000 , 0.802562 )
( 13.700000 , 0.800641 )
( 13.800000 , 0.798704 )
( 13.900000 , 0.796753 )
( 14.000000 , 0.794787 )
( 14.100000 , 0.792805 )
( 14.200000 , 0.790807 )
( 14.300000 , 0.788793 )
( 14.400000 , 0.786762 )
( 14.500000 , 0.784714 )
( 14.600000 , 0.782648 )
( 14.700000 , 0.780565 )
( 14.800000 , 0.778463 )
( 14.900000 , 0.776342 )
( 15.000000 , 0.774201 )
( 15.100000 , 0.772041 )
( 15.200000 , 0.769859 )
( 15.300000 , 0.767657 )
( 15.400000 , 0.765433 )
( 15.500000 , 0.763187 )
( 15.600000 , 0.760917 )
( 15.700000 , 0.758624 )
( 15.800000 , 0.756306 )
( 15.900000 , 0.753962 )
( 16.000000 , 0.751592 )
( 16.100000 , 0.749196 )
( 16.200000 , 0.746771 )
( 16.300000 , 0.744317 )
( 16.400000 , 0.741833 )
( 16.500000 , 0.739318 )
( 16.600000 , 0.736771 )
( 16.700000 , 0.734189 )
( 16.800000 , 0.731573 )
( 16.900000 , 0.728921 )
( 17.000000 , 0.726230 )
( 17.100000 , 0.723500 )
( 17.200000 , 0.720728 )
( 17.300000 , 0.717914 )
( 17.400000 , 0.715053 )
( 17.500000 , 0.712146 )
( 17.600000 , 0.709188 )
( 17.700000 , 0.706177 )
( 17.800000 , 0.703111 )
( 17.900000 , 0.699987 )
( 18.000000 , 0.696800 )
( 18.100000 , 0.693548 )
( 18.200000 , 0.690225 )
( 18.300000 , 0.686828 )
( 18.400000 , 0.683351 )
( 18.500000 , 0.679789 )
( 18.600000 , 0.676134 )
( 18.700000 , 0.672380 )
( 18.800000 , 0.668517 )
( 18.900000 , 0.664537 )
( 19.000000 , 0.660427 )
( 19.100000 , 0.656176 )
( 19.200000 , 0.651766 )
( 19.300000 , 0.647180 )
( 19.400000 , 0.642394 )
( 19.500000 , 0.637381 )
( 19.600000 , 0.632106 )
( 19.700000 , 0.626522 )
( 19.800000 , 0.620571 )
( 19.900000 , 0.614168 )
( 20.000000 , 0.607196 )
};
\addplot [thick, green!80!blue]coordinates{( 0.100000 , 0.991160 )
( 0.200000 , 0.989903 )
( 0.300000 , 0.988644 )
( 0.400000 , 0.987381 )
( 0.500000 , 0.986115 )
( 0.600000 , 0.984846 )
( 0.700000 , 0.983574 )
( 0.800000 , 0.982298 )
( 0.900000 , 0.981019 )
( 1.000000 , 0.979737 )
( 1.100000 , 0.978451 )
( 1.200000 , 0.977162 )
( 1.300000 , 0.975870 )
( 1.400000 , 0.974574 )
( 1.500000 , 0.973275 )
( 1.600000 , 0.971972 )
( 1.700000 , 0.970666 )
( 1.800000 , 0.969356 )
( 1.900000 , 0.968043 )
( 2.000000 , 0.966726 )
( 2.100000 , 0.965405 )
( 2.200000 , 0.964081 )
( 2.300000 , 0.962753 )
( 2.400000 , 0.961421 )
( 2.500000 , 0.960086 )
( 2.600000 , 0.958746 )
( 2.700000 , 0.957403 )
( 2.800000 , 0.956057 )
( 2.900000 , 0.954706 )
( 3.000000 , 0.953351 )
( 3.100000 , 0.951992 )
( 3.200000 , 0.950630 )
( 3.300000 , 0.949263 )
( 3.400000 , 0.947893 )
( 3.500000 , 0.946518 )
( 3.600000 , 0.945139 )
( 3.700000 , 0.943756 )
( 3.800000 , 0.942369 )
( 3.900000 , 0.940977 )
( 4.000000 , 0.939581 )
( 4.100000 , 0.938181 )
( 4.200000 , 0.936777 )
( 4.300000 , 0.935368 )
( 4.400000 , 0.933955 )
( 4.500000 , 0.932537 )
( 4.600000 , 0.931115 )
( 4.700000 , 0.929688 )
( 4.800000 , 0.928256 )
( 4.900000 , 0.926820 )
( 5.000000 , 0.925380 )
( 5.100000 , 0.923934 )
( 5.200000 , 0.922484 )
( 5.300000 , 0.921029 )
( 5.400000 , 0.919568 )
( 5.500000 , 0.918103 )
( 5.600000 , 0.916634 )
( 5.700000 , 0.915159 )
( 5.800000 , 0.913678 )
( 5.900000 , 0.912193 )
( 6.000000 , 0.910703 )
( 6.100000 , 0.909207 )
( 6.200000 , 0.907706 )
( 6.300000 , 0.906200 )
( 6.400000 , 0.904688 )
( 6.500000 , 0.903171 )
( 6.600000 , 0.901648 )
( 6.700000 , 0.900120 )
( 6.800000 , 0.898586 )
( 6.900000 , 0.897046 )
( 7.000000 , 0.895501 )
( 7.100000 , 0.893949 )
( 7.200000 , 0.892392 )
( 7.300000 , 0.890829 )
( 7.400000 , 0.889259 )
( 7.500000 , 0.887684 )
( 7.600000 , 0.886102 )
( 7.700000 , 0.884515 )
( 7.800000 , 0.882920 )
( 7.900000 , 0.881320 )
( 8.000000 , 0.879713 )
( 8.100000 , 0.878099 )
( 8.200000 , 0.876479 )
( 8.300000 , 0.874852 )
( 8.400000 , 0.873218 )
( 8.500000 , 0.871577 )
( 8.600000 , 0.869929 )
( 8.700000 , 0.868274 )
( 8.800000 , 0.866612 )
( 8.900000 , 0.864943 )
( 9.000000 , 0.863267 )
( 9.100000 , 0.861583 )
( 9.200000 , 0.859891 )
( 9.300000 , 0.858192 )
( 9.400000 , 0.856485 )
( 9.500000 , 0.854770 )
( 9.600000 , 0.853047 )
( 9.700000 , 0.851316 )
( 9.800000 , 0.849577 )
( 9.900000 , 0.847829 )
( 10.000000 , 0.846074 )
( 10.100000 , 0.844309 )
( 10.200000 , 0.842536 )
( 10.300000 , 0.840754 )
( 10.400000 , 0.838963 )
( 10.500000 , 0.837164 )
( 10.600000 , 0.835354 )
( 10.700000 , 0.833536 )
( 10.800000 , 0.831708 )
( 10.900000 , 0.829870 )
( 11.000000 , 0.828023 )
( 11.100000 , 0.826166 )
( 11.200000 , 0.824298 )
( 11.300000 , 0.822421 )
( 11.400000 , 0.820532 )
( 11.500000 , 0.818634 )
( 11.600000 , 0.816724 )
( 11.700000 , 0.814804 )
( 11.800000 , 0.812872 )
( 11.900000 , 0.810929 )
( 12.000000 , 0.808974 )
( 12.100000 , 0.807008 )
( 12.200000 , 0.805029 )
( 12.300000 , 0.803039 )
( 12.400000 , 0.801036 )
( 12.500000 , 0.799020 )
( 12.600000 , 0.796991 )
( 12.700000 , 0.794950 )
( 12.800000 , 0.792894 )
( 12.900000 , 0.790826 )
( 13.000000 , 0.788743 )
( 13.100000 , 0.786646 )
( 13.200000 , 0.784534 )
( 13.300000 , 0.782408 )
( 13.400000 , 0.780266 )
( 13.500000 , 0.778110 )
( 13.600000 , 0.775937 )
( 13.700000 , 0.773748 )
( 13.800000 , 0.771543 )
( 13.900000 , 0.769320 )
( 14.000000 , 0.767081 )
( 14.100000 , 0.764824 )
( 14.200000 , 0.762548 )
( 14.300000 , 0.760255 )
( 14.400000 , 0.757942 )
( 14.500000 , 0.755610 )
( 14.600000 , 0.753258 )
( 14.700000 , 0.750885 )
( 14.800000 , 0.748492 )
( 14.900000 , 0.746076 )
( 15.000000 , 0.743639 )
( 15.100000 , 0.741179 )
( 15.200000 , 0.738696 )
( 15.300000 , 0.736189 )
( 15.400000 , 0.733656 )
( 15.500000 , 0.731099 )
( 15.600000 , 0.728515 )
( 15.700000 , 0.725904 )
( 15.800000 , 0.723265 )
( 15.900000 , 0.720597 )
( 16.000000 , 0.717900 )
( 16.100000 , 0.715171 )
( 16.200000 , 0.712411 )
( 16.300000 , 0.709618 )
( 16.400000 , 0.706790 )
( 16.500000 , 0.703927 )
( 16.600000 , 0.701028 )
( 16.700000 , 0.698090 )
( 16.800000 , 0.695112 )
( 16.900000 , 0.692093 )
( 17.000000 , 0.689031 )
( 17.100000 , 0.685924 )
( 17.200000 , 0.682769 )
( 17.300000 , 0.679566 )
( 17.400000 , 0.676311 )
( 17.500000 , 0.673002 )
( 17.600000 , 0.669635 )
( 17.700000 , 0.666210 )
( 17.800000 , 0.662721 )
( 17.900000 , 0.659165 )
( 18.000000 , 0.655539 )
( 18.100000 , 0.651838 )
( 18.200000 , 0.648057 )
( 18.300000 , 0.644191 )
( 18.400000 , 0.640235 )
( 18.500000 , 0.636181 )
( 18.600000 , 0.632022 )
( 18.700000 , 0.627750 )
( 18.800000 , 0.623355 )
( 18.900000 , 0.618826 )
( 19.000000 , 0.614150 )
( 19.100000 , 0.609312 )
( 19.200000 , 0.604294 )
( 19.300000 , 0.599075 )
( 19.400000 , 0.593629 )
( 19.500000 , 0.587924 )
( 19.600000 , 0.581921 )
( 19.700000 , 0.575566 )
( 19.800000 , 0.568792 )
( 19.900000 , 0.561504 )
( 20.000000 , 0.553567 )
};
\end{axis}
\end{tikzpicture}
\caption{Voltages in constant power terminals for different percentage of load}
\label{fig:voltages_cargabilidad}
\end{figure}

The power flow was also evaluated under high-load conditions.  Power in each $k\in\mathcal{P}$ was increased until critical voltage was achieved, in the same manner as in the voltage stability studies for ac power systems. Results are depicted in figs \ref{fig:cargabilidad} and \ref{fig:voltages_cargabilidad}. In normal operative conditions, the successive approximation algorithm converges in less than 5 iterations but the number of iterations increase as the system approaches the critical point. As $P_{max}$ increases, the contraction constant $\alpha$ tends to 1. Nevertheless, convergence is guaranteed in a finite number of iterations even in these extreme conditions. Notice that ac systems require other methodologies such as the continuation power flow for calculations close to maximum load limit. 

\section{Conclusions}

The power flow in LVDC grids was analyzed using the Banach fixed-point theorem. Convergence and uniqueness of the solution was demonstrated under practical considerations.  Simulation results shown that a successive approximations algorithm converges even in operative points close to the voltage collapse.  

\section*{References}

\bibliographystyle{elsarticle-num}
\bibliography{DC_distribution}
%
%\vfill
%\parpic{\includegraphics[width=1.5in,clip,keepaspectratio]{Foto_Alejandro.jpg}}
%\noindent {\bf Alejandro Garc\'es}  was born in Pereira, Colombia, in 1981. He received the Bachelor and Master degree in electrical engineering from Universidad Tecnol\'ogica de Pereira, Colombia in 2006, and the Ph.D. degree from the Norges Teknisk Naturvitenskapelige Universitet (NTNU), Norway, in 2012.  He is currently Associated Professor at the Department of Electric Power Engineering, Universidad Tecnol\'ogica de Pereira. His research interests include renewable energies, mathematical optimization, power system dynamics and HVDC transmission.

\end{document}